\documentclass[12pt]{article}
\usepackage{a4wide,amsmath,amssymb,hyperref}
\usepackage{xcolor,natbib}
\usepackage{verbatim}
\newcommand{\ra}{R_{\alpha}}
\newcommand{\La}{\Lambda}
\newcommand{\De}{\Delta}

\newcommand{\p}{\psi}

\newcommand{\ma}{\beta_\alpha}

\newtheorem{theorem}{Theorem}

\newtheorem{lemma}[theorem]{Lemma}

\newtheorem{thdelcorr}{\rm {\bf Proposition \!\! DEL}\!\!}

\renewcommand{\Bbb}[1]{\mathbb{#1}}
\newcommand{\I}{{\Bbb I}}         
\newcommand{\N}{{\Bbb N}}         
\newcommand{\Q}{{\Bbb Q}}         
\newcommand{\R}{{\Bbb R}}         
\newcommand{\Rp}{(0,+\infty)}    
\newcommand{\Z}{{\Bbb Z}}         

\newcommand{\ds}{\displaystyle}

\def\Ja{Jarn{\'\i}k}
\def\JB{Jarn{\'\i}k-Besicovitch}
\def\K{Khintchine}

\def\Jarnik{Jarn\'\i k}
\def\MTP{Mass Transference Principle}
\def\G{{\rm {\bf G(\alpha)}}}
\newcommand{\cA}{W}

\newcommand{\cC}{{\cal C}}

\newcommand{\cH}{{\cal H}}

\newcommand{\cK}{{\cal K}}
\newcommand{\cL}{{\cal L}}
\newcommand{\cM}{{\cal M}}

\newcommand{\cQ}{{\cal Q}}
\newcommand{\cR}{{\cal R}}

\newcommand{\cV}{{\cal V}}

\newcommand{\vv}[1]{{\mathbf{#1}}}

\def\q{\mathbf{q}}

\def\Bad{{\rm {\bf Bad}}}

\renewcommand{\emph}[1]{{\it#1}}

\begin{document}

\title{The dimension of well approximable numbers}

\author{Victor Beresnevich\footnote{Supported by EPSRC grant EP/Y016769/1}\and
Sanju Velani}

\date{}

\maketitle

\begin{abstract}
In this survey article, we explore a central theme in Diophantine approximation inspired by a celebrated result of Besicovitch  on the Hausdorff dimension of well approximable real numbers. 
We outline some of the key developments stemming from Besicovitch's result, with a focus on the  Mass Transference Principle, Ubiquity and Diophantine approximation on manifolds and fractals.
We highlight the subtle yet profound connections between number theory and fractal geometry, and discuss several open problems at their intersection.
\end{abstract}

{\footnotesize \noindent 2020 Mathematics Subject Classification:\\ Primary: 11J83, 11K60, 11J13, 28A80, 28A78; Secondary: 11K55, 11J71, 42B10, 43A46}

\section{The Jarn\'ik-Besicovitch theorem}\label{sec:Dirichl-Roth-metrica-theory}

In this survey article, we explore a central theme in Diophantine approximation, situated within the broader context of fractal geometry and metric number theory, and inspired by a celebrated result of Besicovitch \cite{MR1574327}, now commonly referred to as the Jarník–Besicovitch theorem. Broadly speaking, given a parameter $\tau > 0$, this theorem determines the fractal dimension of the set of real numbers that are \textit{$\tau$--well approximable};  that is, those real numbers $x$ for which there exist infinitely many $(p, q) \in \mathbb{Z} \times \mathbb{N}$ satisfying
\begin{equation}\label{eq1}
    |x  - p/q|<q^{-\tau-1}   \, . 
\end{equation}
The parameter $\tau$ governs the ``rate'' of  rational approximations $p/q$ to $x$ -- the larger the value of $\tau$, the better the rational approximation one seeks.
It is straightforward to verify that the set of $\tau$--well approximable numbers is invariant under translations by integers. Therefore, without loss of generality, we will restrict our attention (following Besicovitch) to the unit interval $\I := [0,1)$.  With this in mind, we let   
  $$
  W(\tau): =
\{x\in \I \colon \|qx\|<   q^{-\tau} \text{ for infinitely many }
q\in\N \}
$$
denote the set of $\tau$–well approximable numbers in $\I$.  Here and throughout, $\|x\|$ denotes the
distance of $x$ from the nearest integer.

A foundational result in  Diophantine approximation, Dirichlet’s theorem,  states that  given any real number $x$ and any $N \in \N$,  there exists  $(p, q) \in \mathbb{Z} \times \mathbb{N}$ such that 
\begin{equation}\label{Dirn=1}
    \| qx \|   < N^{-1}     \qquad \text{ and }  \qquad \      1 \le q \le N \, . 
\end{equation}
A straightforward corollary of this is that  every real number is approximable by rationals with rate  $\tau = 1$ and thus
\begin{equation}\label{eq2}
 W(1)=\I\,.
 \end{equation}
The Jarn\'ik–Besicovitch theorem reveals the much finer structure of the sets $W(\tau)$ for $\tau > 1$. Although these sets are Lebesgue null, as was already shown by Borel~\cite{MR1511714} in 1912, they are far from negligible in a geometric sense. Their richness is captured by their Hausdorff dimension  -- a key concept in fractal geometry that extends the classical notion of dimension beyond integer values.  Throughout, given a subset $X$ of $\R^n$  we denote by    $ \dim X$  the Hausdorff dimension of  $X$.

\begin{theorem}[\Ja-Besicovitch] \label{JBT} For any $\tau\ge1$, we have that
$$
 \dim W(\tau) = {\frac2{\tau + 1}}\,.
$$
\end{theorem}
This result, established by Jarn\'ik \cite{Jarnik29} in 1929 and independently rediscovered by Besicovitch \cite{MR1574327} in 1934 in a broader geometric context, marked a foundational step in the development of metric Diophantine approximation. It highlights the subtle yet profound connections between number theory and fractal geometry -- connections that continue to inspire a wide range of modern research. In this article we outline some of the key developments stemming from this original discovery and highlight its influence on the wider field. While this article is not intended to be a comprehensive survey of the area encompassing these developments, we hope that it provides sufficient references and context to guide interested readers toward further research and open problems that remain both intricate and compelling.

Before delving deeper into the implications and generalisations of the \JB{} theorem, we begin by introducing the concepts of Hausdorff measure and Hausdorff dimension, both for completeness and to establish the notational framework that will be used throughout the remainder of this article.

\section{Hausdorff measure and dimension} \label{DB}

Let $X$ be  a
subset of $\R^n$. Given a ball $B$ in $\R^n$, let $r(B)$ denote the
radius of $B$. For $\rho > 0$, a countable collection
$ \left\{B_{i} \right\} $ of balls in $\R^n$ with $r(B_i) \le \rho $
for each $i$ such that $F \subset \bigcup_{i} B_{i} $ is called a {\em
  $ \rho $-cover\index{rho-cover} for $F$}.   Let $s$ be a non-negative number and define
$$
{\cal H}^{s}_{\rho} (X) \, := \, \inf \ \sum_{i} r(B_i)^s,
$$
where the infimum is taken over all $\rho$-covers of $X$. The {\it $s$--dimensional
Hausdorff measure}   $ {\cal H}^{s} (X)$ of $X$ is defined by
$$ {\cal H}^{s} (X) := \lim_{ \rho \rightarrow 0} {\cal
H}^{s}_{\rho} (X) \; = \; \sup_{\rho > 0 } {\cal H}^{s}_{\rho} (X)
\;
$$
and the {\it Hausdorff dimension} by
$
\dim X:=  \inf \left\{ s\ge0 : \cH^{s} (X) =0 \right\}   \,.
$
Note that when $s$ is a positive integer, $\cH^s$ is a constant multiple of
Lebesgue measure in $\R^s$.  Indeed, when $s=1$  the two measures were equal if in the definition of ${\cal H}^{s}$  we had used the diameter of a ball instead of its radius. This clearly makes no difference to the definition of dimension.     Further details regarding Hausdorff measure
and dimension can be found in~\cite{MR867284,MR1333890}.

From a historical perspective,  it is well worth highlighting  the following obvious but nevertheless powerful consequence of the \Ja-Besicovitch theorem.    The above  notion of fractal dimension  was introduced by Hausdorff \cite{MR1511917} in 1918,  based on a construction of Carath\'{e}odory. One of the most fundamental and best-known examples of a fractal set is the famous middle third Cantor set, which we will denote by $\cK$. As shown by Hausdorff in the same paper, it has dimension  $\dim \cK = \log2/\log3$. In short, $\cK$ is the limit set obtained by removing the open middle third from the unit interval $[0,1]$, and then recursively repeating this process with each of the remaining subintervals -- ad infinitum.
In arithmetical terms, $\cK$  consists of all real numbers in the unit interval whose base $3$ expansions do not contain the digit $1$ -- that is, only digits $0$ and $2$ appear.
Clearly,  it is possible to generalise the construction of $\cK$ as follows.  Given $R,M\in\N$ such that $M\le R-1$,  partition $[0,1]$  into $R$ equal closed subintervals and remove the interior of any $M$ of them. Then repeat the process with the remaining closed segments - ad infinitum.
The associated limit set  has Hausdorff dimension  $ \log(R-M)/\log R$, see, for example, \cite[Lemma~1.7.3]{MR3618787}.
The upshot of this is that we can construct continuum many ``natural'' Cantor-type sets each with a distinct fractal dimension between zero and one. Of course, not every real number in $(0,1)$ is attained as a dimension via this construction.  The \Ja-Besicovitch theorem shows that the dimension of the set $W(\tau)$  varies continuously with the parameter $\tau$ and thus attains every value in $(0,1]$.   To the best of our knowledge this is the first example of a naturally occurring sets with such a desirable and rich dimensional property.

Continuing with the above ``first example''  theme,  recall that Liouville in 1844 constructed the first explicit  examples of transcendental numbers; i.e., irrational real numbers that are  not the root of a non-zero polynomial with integer (or, equivalently, rational) coefficients.  In terms of the set $W(\tau)$,  the set $\cL$  of Liouville numbers can be characterized as   the set of $\alpha\in\R\setminus\Q$  such that $  \alpha  \in W(\tau) $ for any  $\tau>0$.  Thus, by definition, we have that
$
\mathcal{L}\subset W(\tau)
$
for any $\tau\ge1$, and so, by the monotonicity property of Hausdorff dimension and the  \Ja-Besicovitch theorem, it follows that 
$
0\le \dim\mathcal{L}\le \dim W(\tau)=\frac2{\tau+1}  \, . 
$
Since $\tau$ can be taken arbitrarily large, this means that $\dim\mathcal{L}=0$.  To the best of our knowledge this is the first example of a naturally occurring uncountable set with zero dimension.    Note that in the above argument we only need to use  the  upper bound for   $ \dim W(\tau) $.  We now demonstrate that this  is a simple consequence of the set theoretic  nature of the sets under consideration.

\subsection{The upper bound in the \JB{} theorem}

We start with the important observation that $W(\tau)$ is a $\limsup$
set of the  balls given by \eqref{eq1}  with centres at rational points $p/q\in\I$.  More formally, 
\begin{equation*}
  W(\tau) =  \limsup_{q \to \infty} A_q(\tau) := \textstyle{\bigcap_{t=1}^{\infty} \bigcup_{q=t}^{\infty}} A_q(\tau)\, ,
\end{equation*}
where   for   a fixed $q \in \N$
\begin{eqnarray} \label{slv101}
  A_q (\tau)  &  :=   &    \textstyle{\bigcup_{p=0}^q} B (p/q,
                     q^{-\tau-1} ) \ \cap \I\,    \ :=  \   \{  x \in \I : \|qx\| \, < \,  q^{-\tau} \}   \, .
\end{eqnarray}
Now notice that for each  $t\in \N$, we have that 
$
  W(\tau)    \subset      \textstyle{\bigcup_{q=t}^{\infty}} A_q(\tau) \,
$; 
i.e.\ for each $t$, the collection of balls $B(p/q, q^{-\tau-1}) $
associated with the sets $A_q(\tau): q=t, t+1, \dots $ forms a cover
for $W(\tau)$.
Further, for any $\rho > 0 $, if we choose $t$ large enough so that
$\rho > t^{-\tau-1} $, then  the balls in the collection $\{A_q(\tau)\}_{q \ge t}$
form a $\rho$ cover of $W(\tau)$. Then, 
\begin{equation*}
  \mathcal{H}_{\rho}^s\big(W(\tau)\big)  \ \le  \ \textstyle{\sum_{q = t}^{\infty}}  \ (q+1) q^{-(\tau+1)s}    \to 0
\end{equation*}
as $t \to \infty$ (i.e.\ $\rho \to 0$) if $ s > s_0:=2/(1 + \tau)$.   Hence, it follows from the  definition of Hausdorff measure and dimension that $\mathcal{H}^s\big(W(\tau)\big)   = 0 $  if $ s > s_0$ and  in turn that $ \dim (W(\tau)) \le  s_0$ as desired.   

\medskip

The preceding  argument can be adapted,  with straightforward modifications,  to establish the following result for general $\limsup$ sets.

\begin{lemma}[The Hausdorff-Borel-Cantelli Lemma]\label{HCL}
Let $\{B_k\}$ be a sequence of balls in $\R^n$ with the radii  $r_k$ and let
$$
\textstyle E_\infty=\limsup\limits_{k\to\infty}B_k:=\bigcap_{t=1}^\infty\bigcup_{k\ge t}B_k\,.
$$
Suppose that for some $s>0$ we have that
\begin{equation}\label{eq3}
\sum_{k=1}^\infty r_k^s<\infty\,.
\end{equation}
Then $\cH^s(E_\infty)=0$, in particular, $\dim (E_\infty)\le s$.
\end{lemma}

\noindent
In the case $s=n$, the Hausdorff measure $\cH^s$ is equal to, up to a constant factor, $n$-dimensional Lebesgue outer measure. Thus, within the setting under consideration, Lemma~\ref{HCL} coincides with the standard Borel-Cantelli lemma from probability theory.
To the best of our knowledge, Lemma~\ref{HCL} is the key tool for determining upper bounds on the Hausdorff dimension of a wide class of limsup sets that naturally arise in number theory. The primary challenge in applying this lemma is to find an efficient representation of the limsup sets under consideration that ensures the convergence of \eqref{eq3}.

Lemma~\ref{HCL} was coined `the Hausdorff-Cantelli Lemma' by Bernik and Dodson \cite{MR1727177}. 
The terminology appears to serve as an acknowledgment of the ideas involved in the proof --  it does not appear that either Hausdorff or Cantelli did any formal work to adapt the Borel-Cantelli Lemma to Hausdorff measures. Henceforth, if any name is to be used, it would likely be more appropriate to refer to the lemma as the Hausdorff-Borel-Cantelli Lemma. We also note that the Borel–Cantelli type argument, which underpins  the proof of Lemma~\ref{HCL}, has been used implicitly for decades by multiple authors without being formalised into a self-contained statement. In particular, it can be recovered from Besicovitch's proof of statement (i) in \cite[p.127]{MR1574327}.

We now turn our attention to discussing the lower bound for  $\dim W(\tau)$, which constitutes the main substance of the \JB{} theorem. Our discussion will be framed within the context of general $\limsup$  sets -- yet another illustration of the theorem’s far-reaching influence.
In what follows, we present two general approaches for establishing lower bounds for $\limsup$ sets -- Mass Transference and Ubiquity, both of which have been influenced by the \JB{} theorem and its various generalisations.

\section{Lower bounds via Mass Transference}

In this section, we present the Mass Transference Principle, introduced in \cite{MR2259250}, and demonstrate its applications in the context of the \JB{} theorem. For the sake of clarity,  rather than working in arbitrary metric spaces, it is appropriate to restrict our attention to $\R^n$ as in the original formulation  appearing  in the main body of \cite{MR2259250}.  Throughout,  $B(x,r)   \subset \R^n $ will denote a ball with centre $x$ and radius $r$.

\subsection{Mass Transference Principle}\label{mtp0}

Given $s>0$, define the following transformation
on balls in $\R^n$:
\begin{equation}\label{B^s}
\textstyle B=B(x,r)\mapsto B^s:=B(x,r^{s/n}) \ .
\end{equation}
Clearly $B^n=B$. The following statement
is the key to obtaining lower bounds for Hausdorff measure and dimension of limsup sets of balls. It is very much influenced by the desire to obtain an analogue of the \JB{} theorem for general limsup sets.

\begin{theorem}[\MTP{} \cite{MR2259250}]\label{MTPthm}
Let $\{B_k\}_{k\in\N}$ be a sequence of balls in $\R^n$ with
radii $r_k\to 0$ as $k\to\infty$. Let $s>0$ and 
$\Omega$ be any bounded open subset in $\R^n$. Then
\begin{equation}\label{MTP0}
\cH^n\big(\/\Omega\cap\limsup_{k\to\infty}B^s_k{}\,\big)=\cH^n(\Omega) \quad\Longrightarrow\quad
\cH^s\big(\/\Omega\cap\limsup_{k\to\infty}B^n_k\,\big)=\cH^s(\Omega) \ .
\end{equation}
\end{theorem}

Note that for $s>n$ the statement is trivial since $\cH^s(\R^n)=0$. Thus the main substance of Theorem~\ref{MTPthm} is the case $s<n$. Recall that $\cH^n$ is comparable to the
$n$-dimensional Lebesgue measure and so the \MTP{} allows us to transfer Lebesgue measure statements to Hausdorff measure statements. In relation to the \JB{} theorem, the lower bound easily follows from \eqref{eq2} via the \MTP, which we apply to $\Omega=(0,1)$ and the sequence of balls $B(p/q, q^{-\tau-1}) $
associated with the sets $A_q(\tau)$ defined by \eqref{slv101}. Indeed, on taking $s=\frac{2}{\tau+1}$ and using  definition of $B^s$ -- see \eqref{B^s}, we observe that $B^s(p/q, q^{-\tau-1})=B(p/q, q^{-2})$. Therefore, since $p$ and $q$ are not necessarily coprime,
we get that $\limsup B^s(p/q, q^{-\tau-1})=W(1)$, which, by \eqref{eq2}, is of full Lebesgue measure in $\Omega=(0,1)$. This provides the hypothesis of the Mass Transference Principle, that is the left hand side of \eqref{MTP0}. Furthermore, as we have seen already seen, $\limsup B(p/q, q^{-\tau-1})=\limsup A_q(\tau)=W(\tau)$. Therefore, by \eqref{MTP0}, we obtain that $\cH^s(W(\tau))=\cH^s(\Omega\cap W(\tau))=\cH^s(\Omega)>0$ for $s\le 1$ (equivalently $\tau\ge 1$). Furthermore, if $\tau>1$ then $s<1$, and we get that  $\cH^s(\Omega) =\infty$. Thus, the Hausdorff measure at the critical exponent is infinite, that is 
$$\cH^s(W(\tau))=\infty\quad\text{for $s=\tfrac2{\tau+1}   \ \ (\tau > 1) $}\,.
$$

\subsection{General approximation functions}

In the above section we have been dealing with $\tau$-approximable numbers associated with
the error of
approximation given by
$q^{-\tau}$. More generally,
given a function  $\psi:\N\to[0,+\infty)$, we consider the set of $\p$-well approximable numbers
  $$
  W(\psi): =
\{x\in \I \colon \|qx\|<\psi(q) \text{ for infinitely many }
q\in\N \}  \ .
$$
Regarding this general sets we have the following elegant criterion concerning its size. 

\begin{theorem}[The Khintchine-Jarn\'{\i}k theorem]
\label{khijarmeas} Let $s>0$. Then
$$
\cH^s\left(W(\psi)\right)=\left\{\begin{array}{cl} 0 &
\ds\text{if } \;\;\; \sum_{r=1}^{\infty}  \ \; r \,
\left(\frac{\psi(r)}r\right)^s
 <\infty \; ,\\[3ex]
\cH^s(\I) & \ds\text{if } \;\;\; \sum_{r=1}^{\infty} \  \; r
\, \left(\frac{\psi(r)}r\right)^s =\infty   \  \text{ and $\psi$
is monotonic}.
\end{array}\right.$$
\end{theorem}

{For $s>1$ the statement is trivial since both sides vanish.} The $s=1$ case of this statement was discovered by  \K~\cite{MR1512207} and the $s<1$ was proved by Jarn\'ik \cite{Jarnik29}, albeit both had additional constrains on $\psi$.  In view of the Mass Transference Principle it is easily verified  (see \cite{MR2259250} or \cite[\S3.4]{MR3618787})   that
\begin{center} Khintchine's Theorem $ \hspace{4mm} \Longrightarrow
\hspace{4mm} $ Jarn\'{\i}k's Theorem.  \end{center} 
Thus, the
Lebesgue theory of $W(\psi)$ underpins the general Hausdorff
theory.

The original discovery of the Mass Transference Principle was very much influenced by the desire to deal with more sophisticated limsup sets. In particular, the main motivation was to establish the Duffin-Schaeffer conjecture \cite{MR4859,MR4125453} and its higher dimensional generalisation \cite{MR1099767} for Hausdorff measures. 
Indeed, in 1941 Duffin and Schaeffer \cite{MR4859} constructed an example of a
non-monotonic function $\psi$ for which the sum
$\sum_q \psi(q)$ diverges but $W(\psi)$ is null. This examples shows that the monotonicity condition is essential in the statement of Khintchine's theorem.   In the same paper they 
provided an alternative  statement for arbitrary $\psi$.  This soon  became known as the  Duffin-Schaeffer conjecture and involved considering the alternative set
  $$
  W'(\psi): =
\{x\in \I \colon |qx-p|<\psi(q) \text{ for infinitely many }
\text{coprime pairs }(q,p)\}  \ .
$$
The coprimality condition imposed on $p$ and $q$ plays a crucial role in formulating the Duffin–Schaeffer conjecture by ensuring that the approximating rational fractions $p/q$ are reduced. In particular, this means that the sequence of rational numbers $p/q$ generated by the pairs $(p,q)$ does not have repetitions, which are otherwise exploited in the Duffin-Schaeffer counterexample. We also note that for each $q\in\N$ there are exactly $\varphi(q)$ rational numbers $p/q$ written as reduced fractions and lying in the unit interval $\I$. Here and elsewhere $\varphi$ is Euler's totient function.

Regarding the sets $W'(\psi)$ we have the following criterion that generalises Theorem~\ref{khijarmeas} to arbitrary (not necessarily monotonic) functions $\psi$:

\begin{theorem}[The Duffin-Schaeffer conjecture for Hausdorff measures]
\label{GDS} Let $s>0$. Then
$$
\cH^s\left(W'(\psi)\right)=\left\{\begin{array}{cl} 0 &
\ds\text{if } \;\;\; \sum_{q=1}^{\infty}  \ \; \varphi(q) \,
\left(\frac{\psi(q)}q\right)^s
 <\infty \; ,\\[3ex]
\cH^s(\I) & \ds\text{if } \;\;\; \sum_{q=1}^{\infty} \  \; \varphi(q)
\, \left(\frac{\psi(q)}q\right)^s =\infty   \ .
\end{array}\right.$$
\end{theorem}

{Similarly to Theorem~\ref{khijarmeas}, for $s>1$ the statement is trivial.} The $s=1$ case of this statement, which corresponds to Lebesgue measure, is the content of the Duffin-Schaeffer conjecture \cite{MR4859}, originally formulated in 1941.
Until recently, it remained one of the major open problems in metric number theory, see \cite{MR2136100, MR1672558}. A complete proof was finally established in a major breakthrough by Koukoulopoulos and Maynard \cite{MR4125453} in 2020. The case $s<1$ was previously established in \cite{MR2259250}, conditional on the $s=1$ case, through the introduction of the Mass Transference Principle. Consequently, in light of the result of Koukoulopoulos and Maynard, Theorem~\ref{GDS} now holds unconditionally for all $s>0$.

{We note that the setup of Theorem~\ref{GDS} allows to take $\psi$ to be supported on an infinite subset $\cQ\subset\N$ that is $\psi(q)=0$ if $q\not\in\cQ$. In particular, by taking $\psi_{\cQ,\tau}(q)=q^{-\tau}$ for $q\in\cQ$ and $0$ otherwise, and using the well known fact that $q/\log(q)\le \varphi(q)\le q$ for all sufficiently large $q$, we can recover from Theorem~\ref{GDS} the following result originally established by Borosh and Fraenkel \cite{MR308060} in 1972:
\begin{equation}\label{BF}
\dim W(\psi_{\cQ,\tau})=\frac{\nu+1}{\tau+1}\qquad\text{for all }\tau\ge\nu\,,
\end{equation}
where
$$
\nu=\nu_\cQ:=\inf\left\{\lambda>0:\sum_{q\in\cQ}q^{-\lambda}<\infty\right\}\,.
$$
For instance, if $\cQ=\{2^k:k\in\N\}$, then $\nu=0$ and $\dim W(\psi_{\cQ,\tau})=\frac{1}{\tau+1}$ for all $\tau>0$.}

\medskip

There are numerous generalisations of the Mass Transference Principle, including extensions to higher-dimensional real spaces and abstract metric spaces, to more general Hausdorff measures, and to settings involving weighted Diophantine approximation, where the approximating sets are rectangles. Further developments include generalisations to approximating sets of arbitrary shapes, as well as applications to multifractal analysis and inhomogeneous Diophantine approximation, to name a few. In each case, a corresponding version of the \JB{} theorem can be formulated. In what follows we explore higher dimensional generalisations, which will also be useful for a subsequent discussion of Diophantine approximation on manifolds.

\subsection{Higher dimensions: systems of linear forms}

The theme of one-dimensional Diophantine approximation discussed so far admits a natural generalisation to higher dimensions by considering systems of simultaneously small homogeneous linear inequalities.

Let $X=(x_{ij})$ be a real $n
\times m$ matrix, regarded as a point in $\R^{nm}$. Thus,  $X$ is
said to be $\psi$-well approximable if the inequality
$$
  \| \q X \| < \psi(|\q|)
$$
\noindent is satisfied for infinitely many $\q=(q_1,\dots,q_n) \in \Z^n$. Here $\q
X$ is the system
$$
q_1x_{1j} + \dots + q_n x_{n,j} \hspace*{6ex}  ( 1\le j \le m )
$$
of $m$ real linear forms in $n$ variables
and  $\| \q X \| := \max_{1\le j \le m } \| \q X^{(j)}
\|$, where $ X^{(j)}$ is the $j$'th column vector of $X$. Since the
set of $\psi$-approximable points is translation invariant under
integer vectors, we can restrict our attention to the $nm$-dimensional
unit cube $\I^{nm}$. The set of $\psi$-approximable points in
$\I^{nm}$ will be denoted by
$$
\cA_{n,m}(\psi) := \{X\in \I^{nm}:\|\q X \| < \psi(|\q|) {\text {
for infinitely many }} \q \in \Z^n \}\,,
$$
{where $|\vv q|=\max\{|q_1|,\dots,|q_n|\}$.} Furthermore, if $\psi(q)=q^{-\tau}$ for some $\tau>0$, we write $\cA_{n,m}(\tau)$ for $\cA_{n,m}(\psi)$.
Note that when $n=m=1$,  $\cA_{n,m}(\psi)$ and $\cA_{n,m}(\tau)$ reduce to the sets $\cA(\psi)$ and $\cA(\tau)$ respectively, introduced earlier in the one-dimensional setting. 
Also it is worth mentioning the two special cases of the general framework of systems of linear forms: the simultaneous and dual forms of approximation corresponding to $n=1$ and $m=1$ respectively \cite{MR1727177}. In the simultaneous case, which corresponds to approximation by rational points {with the same denominator}, we will use the notation 
\begin{equation}\label{W}
\cA_m(\psi):=\cA_{1,m}(\psi)\qquad\text{and}\qquad\cA_m(\tau):=\cA_{1,m}(\tau)\,,
\end{equation} 
while in the dual case, we will write
\begin{equation}\label{W*}
\cA^*_n(\psi):=\cA_{n,1}(\psi)\qquad\text{and}\qquad\cA^*_n(\tau):=\cA_{n,1}(\tau)\,.
\end{equation}

A generalisation of the \JB{} theorem, which provides the Hausdorff dimension of $\cA_{n,m}(\tau)$, was established by Bovey and Dodson \cite{MR847294}. Specifically, they proved that
\begin{equation}\label{BovDod}
\dim\cA_{n,m}(\tau)=(n-1)m+\frac{n+m}{\tau+1}\qquad\text{for any}\quad \tau\ge n/m\ .
\end{equation}
{Indeed, $\cA_{n,m}(n/m) = \I^{nm}$.} Later, Dodson \cite{MR1184759} obtained the Hausdorff dimension of $\cA_{n,m}(\psi)$ for arbitrary monotonic function $\psi$.
In short, for these sets the quantity $\tau$ in the right hand side of \eqref{BovDod} is replaced by the
lower order at infinity of the function $1/\psi$  
$$
\lambda_{1/\psi} := \liminf_{q \to \infty} \frac{\log 1/\psi(q)}{ \log q} \, .   $$

A far reaching generalisation of the Mass Transference Principle to systems of linear forms was obtained in \cite{MR3798593, MR2264714}. Together with the most general version of the Khintchine-Groshev theorem established in \cite{MR2576284} this gives the following general statement for higher dimensions:

\begin{theorem}
\label{Groshevthm}  Let $\psi:\N\to[0,+\infty)$ and $nm>1$. Then for any $s>0$
$$
    \cH^s(\cA_{n,m}(\psi)) = \begin{cases} 0
      & \ \ds\text{if } \  \sum_{q =1}^\infty q^{m+n-1}\left(\frac{\psi(q)}{q}\right)^{s-m(n-1)}<\infty, \\[3ex]
      \cH^s(\I^{nm})
      & \ \ds\text{if }  \  \sum_{q =1}^\infty q^{m+n-1}\left(\frac{\psi(q)}{q}\right)^{s-m(n-1)}=\infty \ .
                  \end{cases}
$$
\end{theorem}

Due to the counterexample of Duffin and Schaeffer \cite{MR4859} the above theorem is false for $m=n=1$ unless additional conditions are imposed, such as the monotonicity of $\psi$. Thus the result in this form is best possible. For monotonic $\psi$ Theorem~\ref{Groshevthm} was originally established in \cite{MR1468922}, which in turn was a generalisation of Jarn\'ik's theorem \cite{MR1545226} for $n=1$.
We finally note that a version of Theorem~\ref{Groshevthm} for more general Hausdorff measures was also established in \cite{MR3798593,MR2264714}. For further insights and developments in the topic we refer the reader to \cite{MR4909148, MR3989813, MR2508636,  MR2370280, MR3618787, MR2072751,  MR4902008, MR4891430,  MR4776193, MR1760086} and references within.

\section{Ubiquity and the manifolds theory}\label{RSU}

In this section we discuss the so-called ubiquity approach for establishing lower bounds for the Hausdorff measure and dimension of $\limsup$ sets. The concept of a ubiquitous system was originally introduced by Dodson, Rynne and Vickers \cite{MR1067887} in 1990 with the aim of obtaining lower bounds for the Hausdorff dimension of $\psi$-well approximable points lying on manifolds -- see \cite{MR985691} for the first application of this kind. 
The framework builds upon the notion of regular systems introduced by Baker and Schmidt \cite{MR271033} in 1970. In fact, Rynne \cite{MR1203279} demonstrated that every regular system is a ubiquitous system, and that the converse is not always true. We also refer the reader to the survey \cite{MR1975457} on the topic. The notion of ubiquity was later developed in \cite{MR2184760} into a general framework allowing for more refined computations of the Hausdorff measure of $\limsup$ sets and encompassing a wide variety of settings.

We now give a very much simplified version of the definition of local ubiquity from \cite{MR2184760} by restricting ourselves to the ubiquitous systems of points  in $\R^n$ with respect to Lebesgue measure only. This will be enough to highlight the key ideas.

\medskip

\noindent\textbf{Definition (Ubiquitous systems)\,:} 
Let $\Omega$ be an open subset of $\R^n$, $\cR=(\ra)_{\alpha \in J }$ be a sequence of points $\ra$ in $\Omega$ indexed by a countable set $J$, and let $\beta: J \to \Rp : \alpha
\mapsto \ma $ a positive function on $J$ that attaches a `weight' $\ma$  to $\ra$, such that $\#\{\alpha\in J: \ma\le b\}<\infty$ for any $b>0$.
Furthermore, let $k>1$ and $\rho : \Rp \to \Rp $ be a  function with $\rho(r) \to 0 $ as $r
\to \infty $. Define
\begin{equation}\label{Delta}
\De(\rho,k;t) := \bigcup_{\alpha\in J:\ \ma\le k^t} B(\ra,\rho(k^t))
\ .
\end{equation}
Then, the pair $(\cR,\beta)$ is said to be a {\em (local)
ubiquitous system relative to $(\rho,k)$} if for any ball $B$ centred in $\Omega$ of sufficiently small radius one has that
\begin{equation}\label{ub}
  \lambda_n\left( B \cap \De(\rho,k;t) \right) \ \ge \ \kappa \ \lambda_n(B) \quad \mbox{for $t\in\N$ large enough}
\end{equation}
for some constant $\kappa>0$ independent from $B$, where $\lambda_n$ denotes Lebesgue measure on $\R^n$.

\medskip

The ubiquity hypothesis \eqref{ub} formalises a covering principle and hence imposes certain condition on the distribution of the points $\ra$. In fact, as part of its proof, the Mass Transference, that we have discussed earlier, also uses a covering property (Lemma~5 in \cite{MR2259250}) which is equivalent to the left hand side of \eqref{MTP0}. However, unlike ubiquity, the latter provides no control over the size of the balls in the cover. 
Indeed, the additional control over the size of the balls in the ubiquity framework -- the balls in \eqref{Delta} are actually of the same size -- enables us to establish  an analogue of a Khintchine type theorem as well as a  \JB{} type theorem for the $\limsup$ set 
$$\La(\Psi)=\{x\in\Omega:x\in B(\ra,\Psi(\ma))\ \mbox{for\
infinitely\ many\ }\alpha\in J\} \ . $$
associated with the ubiquity framework.  Here $\Psi : \Rp \to \Rp $ is a decreasing function and  following is a simplified version of the key  measure theoretic statements regarding the size of $ \La(\Psi)$ obtained in \cite{MR2184760}.

\begin{lemma}[Ubiquity Lemma]\label{UL}
Suppose that $(\cR,\beta)$ is a
ubiquitous system relative
 to $(\rho,k)$, $\Psi:(0,+\infty)\to(0+\infty)$ is monotonic, and $s>0$. Furthermore, suppose that  $\rho$ satisfies the following
$k$-regularity condition:
there exists a strictly positive constant $\lambda < 1$ such that
for $t$ sufficiently large
\begin{equation} \label{reg0}
\rho(k^{t+1}) \leq \lambda \, \rho(k^t) \ .
\end{equation}
Then, whenever we have the following divergence sum condition
\begin{equation}  \label{unif}
  \sum_{t=1}^{\infty} \
 \frac{\Psi(k^t)^s}{\rho(k^t)^{n} }  \ = \ \infty
\end{equation}
we have that
$${\cal H}^s  \left( \Lambda(\Psi) \right) \ = \ {\cal H}^s
\left( \Omega \right) \ \ .
$$
\end{lemma}

We would like to re-emphasize that the Ubiquity Lemma provides statements concerning both the ambient measure (Lebesgue measure on $\R^n$) and the Hausdorff measures. Observe that within the framework of Mass Transference Principle, the ambient measure case serves as the initial hypothesis -- see \eqref{MTP0}.

We now give a basic example of a ubiquitous system in $\R^n$. Set $\cR=(\vv p/q)_{(\vv p,q)\in J}$, indexed by the subset $J\subset\Z^n\times\N$ such that $\vv p/q\in(0,1)^n$ whenever $(\vv p,q)\in J$, and define $\beta_{\vv p,q}=q$. Then there is $k>1$ such that the pair $(\cR,\beta)$ is a
ubiquitous system in $(0,1)^n$ relative to $(\rho,k)$ with  $\rho(r)=
 k\times r^{-(n+1)/n}$. This follows rather easily from Dirichlet's theorem -- this is 
statement \eqref{Dirn=1} for $n=1$. Indeed, \eqref{Dirn=1} already provides a covering of $(0,1)$ by balls $B(p/q,1/(qN))$ with $q\le N:=k^t$. To get \eqref{ub}, one only needs to show that the balls with $q\le k^{t-1}=N/k$ cover only a relatively small {measure within every open interval $B\subset (0,1)$. Indeed, the measure in question is bounded by
$$
\sum_{q\le N/k}\ \sum_{\substack{0\le p\le q\\[0.5ex] B(p/q,1/(qN))\cap I\neq\varnothing}} \frac{2}{qN}\le \sum_{q\le N/k} \frac{2q|B|+4}{qN}\le \frac{2|B|}{k}+\frac{4}{N}\log\left(\frac Nk+1\right)\le \frac12|B|
$$
for a fixed $k>4$ and all sufficiently large $N$.
As a results, for all large $t$
$$
B\cap \bigcup_{k^{t-1}<q\le k^t} \; \bigcup_{p=0}^q B\left(\frac pq,\frac{k}{k^{2t}}\right)\ \supset\  
B\cap \bigcup_{k^{t-1}<q\le k^t} \; \bigcup_{p=0}^q B\left(\frac pq,\frac{1}{qN}\right)
$$
has measure $\ge\frac12|B|$, and this consequently verifies \eqref{ub} for all large $t$ with $\kappa=1/(2k)$.}

For early non-trivial examples of ubiquitous systems, albeit stated in terms of regular systems, we refer to \cite{MR271033, MR1709049}. These deal with algebraic numbers of a fixed degree $n$. In particular, in \cite{MR1709049} it was shown that real algebraic numbers $\alpha$ of degree $n$ with $\ma:=H(\alpha)^{n+1}$ form a ubiquitous system in any finite interval $\Omega\subset\R$, where $H(\alpha)$ is the naive height of the minimal polynomial of $\alpha$ over $\Z$, {that is, the maximum of the moduli of its coefficients}. Consequently, by using Lemma~\ref{HCL} together with a straightforward counting argument, and applying Lemma~\ref{UL}, we have the following general result regarding the $\limsup$ set
$$
A_{n}(\Psi) := \{x\in \R:|x-\alpha| < \Psi(H(\alpha)) {\text {
for infinitely many }} \alpha\in\R,\; \deg_\Q\alpha=n\}\,.
$$

\begin{theorem}[The Khintchine-Jarn\'{\i}k theorem for algebraic numbers]
\label{Algebraic}  Let $s>0$ and $\Psi:\N\to[0,+\infty)$. Then for any interval $I\subset\R$
$$
    \cH^s(A_{n}(\Psi)\cap I) = \begin{cases} 0
      & \ \ds\text{if } \  \sum_{h =1}^\infty h^{n}\Psi(h)^{s}<\infty, \\[3ex]
      \cH^s(I)
      & \ \ds\text{if }  \  \sum_{h =1}^\infty h^{n}\Psi(h)^{s}=\infty \quad\text{and $\Psi$ is monotonic}\ .
                  \end{cases}
$$
\end{theorem}

For $n=1$, this is essentially Theorem~\ref{khijarmeas}. The ambient measure case $s=1$ was established in \cite{MR1709049} and the Hausdorff measures case $s<1$ was established in \cite{MR2184760}. We also note that Theorem~\ref{Algebraic} implies the following fundamental result on the dimension of $A_n(\tau):=A_n(\Psi_\tau)$, where $\Psi_\tau(h):=h^{-\tau-1}$:
$$
\dim A_n(\tau):=\frac{n+1}{\tau+1}\qquad\text{for any }\tau\ge n\,.
$$
This was originally established by Baker and Schmidt \cite{MR271033} and represents a generalisation of the \JB{} theorem for approximations by algebraic numbers of bounded degree. Moreover, in \cite{MR271033}, Baker and Schmidt established the lower bound in the following fundamental result in transcendental numbers theory: 
\begin{equation}\label{BS}
\dim W^*_n(\tau)\cap\cV_n= \frac{n+1}{\tau+1}\qquad\text{for any }\tau\ge n\,,
\end{equation}
where $\cV_n:=\{(x,\dots,x^n):x\in\R\}$ is a curve in $\R^n$, commonly referred to as a Veronese curve. Here $W^*_n(\tau)$ is given by \eqref{W*}. The upper bound in \eqref{BS} was later established by Bernik \cite{MR729734}. The two papers \cite{MR271033, MR729734} essentially laid the foundations for the use of the Hausdorff dimension in the theory of Diophantine approximation on manifolds which we now proceed to discuss.

\subsection{Diophantine approximation on manifolds}

We first start with the dual theory for manifolds, which concerns results on the Hausdorff dimension of the intersection of $W^*_n(\tau)$ with submanifolds  of $\R^n$.  
A very general result in this direction was established by Dickinson and Dodson \cite{MR1738177} who proved that for every extremal manifold $\cM$ in $\R^n$ the following lower bound holds:
\begin{equation}\label{DD}
\dim W^*_n(\tau)\cap\cM \ge \frac{n+1}{\tau+1}+\dim\cM-1\qquad\text{for any }\tau>n\,.
\end{equation}
Recall that a $d$-dimensional manifold $\cM$ in $\R^n$ is extremal if $\cM\cap W^*_n(\tau)$ is null, that is $\cH^d(\cM\cap W^*_n(\tau))=0$, for any $\tau>n$.  On imposing a stronger condition, namely that of  nondegeneracy, \eqref{DD} was extended to a Khintchine-Jarn\'ik type theorem for divergence \cite{MR2184760}. Recall that, in the analytic case, a  submanifold $\cM$ of $\R^n$ is nondegenerate if it can be locally parametrised by real analytic functions $f_1,\dots,f_n$ such that $1,f_1,\dots,f_n$ are linearly independent over $\R$. For a more general definition of nondegeneracy, see \cite{MR1652916}. 

Going back to \eqref{DD}, Baker \cite{MR485713} proved that for nondegenerate curves in $\R^2$ the bound in \eqref{DD} is exact. Due to \eqref{BS}, it is also exact for $\cM=\cV_n$ -- the Veronese curves. However, beyond $n=2$, the Veronese curves and their linear transformations over $\Q$, not much is known about the complementary upper bounds. Indeed, the following intricate problem remains wide open:

\medskip

\emph{Problem $:$ \ \ Prove that for any analytic nondegenerate submanifold $\cM$ of $\R^n$ \eqref{DD} is sharp, that is}
\begin{equation}\label{DD+}
\dim W^*_n(\tau)\cap\cM = \frac{n+1}{\tau+1}+\dim\cM-1\qquad\text{for any }\tau>n\,.
\end{equation}

\medskip

To date \eqref{DD+} has only ever been verified for a small range of $\tau$, namely $n<\tau<n+c/n$  where  $c\approx 0.25$  \cite{MR2069553}, and for some special classes of manifolds \cite{MR985691} which exclude curves. For instance, we do not have an answer to this problem even in the polynomial case, specifically for manifolds parametrised by polynomials over $\Z$, e.g. for $\cM\subset\R^3$ given by $(x,x^2,x^4)$.

\subsection{Simultaneous approximations on manifolds}

Now we discuss the theory of simultaneous Diophantine approximation on manifolds, that is results concerning the Hausdorff dimension of the intersection of $W_n(\tau)$ with submanifolds of $\R^n$. Here $W_n(\tau)$ is the set introduced in \eqref{W}. 
In general, the dimension of the intersection $W_n(\tau) \cap \cM$ cannot be captured by a single expression -- or even several expressions -- that holds for all non-degenerate manifolds of a fixed dimension and all values of $\tau$. Indeed, for larger values of $\tau$, the dimension is known to depend not only on $\tau$, $n$, and $\dim \cM$, but also on finer geometric and arithmetic properties of $\cM$; see, for example, \cite{MR2874641}. However, when $\tau$ is sufficiently close to the (Dirichlet) exponent $1/n$, we expect a uniform formula depending on $\tau$, $n$ and $\dim\cM$ to hold for most, if not all, non-degenerate manifolds. This leads to the following problem, based on heuristic reasoning involving the expected number of rational points of bounded height near a manifold:

\medskip

\emph{Problem $:$ \ \ Let $1\le d< n$ be integers and $\cM$ be a submanifold in $\R^n$ of dimension $d$ and codimension $m=n-d$. Determine the largest $\tau_{n,d}(\cM)>1/n$ such that} 
\begin{equation}\label{DDSim}
\dim W_n(\tau)\cap\cM = \frac{n+1}{\tau+1}-m\qquad\text{for any }1/n< \tau<\tau_{n,d}(\cM)\,.
\end{equation}

\medskip

This problem can also be considered separately in relation to the corresponding lower and upper bound in \eqref{DDSim}. More precisely, it is of interest to determine the largest $\tau_{n,d}^l(\cM)$ (resp. $\tau_{n,d}^u(\cM)$) such that  \eqref{DDSim} holds for all $1/n<\tau<\tau_{n,d}^l(\cM)$ (resp. $1/n<\tau<\tau_{n,d}^u(\cM)$) with `$\ge$' (resp. `$\le$') instead of `$=$'. It is also of interest to find bounds for $\tau$ that guarantee \eqref{DDSim} (or the upper/lower bound) within a class of manifolds. With this in mind, we define
$$
\tau_{n,d}:=\inf\tau_{n,d}(\cM)\,,\qquad \tau_{n,d}^l:=\inf\tau_{n,d}^l(\cM)\,,\qquad \tau_{n,d}^u:=\inf\tau_{n,d}^u(\cM)\,,
$$
where the infimum is taken over analytic nondegenerate manifolds $\cM$ in $\R^n$ of dimension $d$.
Various results and heuristics lead to the following conjecture, which appears, for example, as Conjectures 2.4 and 2.7 in \cite{MR4652409}:

\medskip

\emph{Conjecture $:$ \ \ For $n>d\ge 1$
$$
\tau_{n,d}=\left\{\begin{array}{cl}
\frac{1}{n-d}\,,&\quad\text{if $d>1$\,,}\\[2ex]
\frac{3}{2n-1}\,,&\quad\text{if $d=1$\,.}
\end{array}\right.
$$}

\medskip

It is known that $\tau_{n,d}^u\le 1/(n-d)$ when $d>1$, that $\tau_{n,d}^l\ge 1/(n-d)$ and that for curves we have that  $\tau_{n,1}^l\ge 3/(2n-1)$ \cite{MR2874641}. The situation was completely resolved for planar curves where $\tau_{2,1}^l=\tau_{2,1}^u=1$, see \cite{MR2373145, MR2242634}. Also, for any hypersurface (that is $d=n-1$) $\cM$ in $\R^n$ with non-zero Gaussian curvature we have that $\tau_{n,n-1}^l(\cM)=\tau_{n,n-1}^u(\cM)=1$, see \cite{MR4132580}. In general, it was proved in \cite{MR4652409} that $\tau_{n,d}^u\ge 1/n+O(n^{-3})$ with an explicitly computable constant. This constant was recently slightly improved in \cite{BerDat,SST}. Also, in the case of the Veronese curves $\cV_n$ Badziahin \cite{MR4881431} established that $\tau_{n,1}(\cV_n)\ge 1/n+0.5n^{-2}+O(n^{-3})$ and most impressively proved that $\tau_{3,1}(\cV_3)\ge 3/5$. The latter, for the first time, attains the values of $\tau$ suggested in the above conjecture for a nondegenerate manifold of codimension other than $1$. For further insights into the topic we refer the reader to recent publications \cite{MR1816807, MR4397037, MR4706444, MR3430242, MR4916720} and references within these and aforementioned papers.

\medskip

{We note that \eqref{DDSim} as well as \eqref{DD+} are consistent with the identity 
\begin{equation}\label{generic}
\operatorname{codim}(A\cap B)=\operatorname{codim}(A)+\operatorname{codim}(B)
\end{equation}
which holds generically, that is to say that the sets in question are ``independent'' -- equivalently ``random''. Determining 
$\dim W_n(\tau)\cap\cM$ for $\tau\ge \tau_{n,d}(\cM)$  seems to depend on the `internal arithmetic' of $\cM$, that is on the number and distribution of rational points lying on or very near $\cM$. Little is known about this problem, however, when $\cM\subset\R^2$ is either the unit circle $x^2+y^2=1$ or the parabola $y=x^2$, or a rational transformation of the circle or the parabola, then we have the following complete answer \cite{MR2373145, MR2604984, MR1816807,MR549277}:
\begin{equation}\label{PC}
\dim(W_2(\tau)\cap \cC)=\max\left\{\;\frac{3}{\tau+1}-m,\ \ \frac{d}{\tau+1}\;\right\}\qquad\text{for all }\tau\ge1/2\,,
\end{equation}
where $m=\operatorname{codim}(\cM)=1$ and $d=\dim(\cM)=1$. In higher dimensions, it would be interesting to determine $\dim(W_n(\tau)\cap\cV_n)$ for all $\tau\ge1/n$, for every Veronese curve $\cV_n$. For further insights and developments regarding this problem, we refer the reader to \cite{MR4881431,MR4089038, MR2791654} and references within.}

\medskip

We finally note that, of course, there is nothing stopping us from replacing manifolds $\cM$ in the above discussion by more  general objects, such as fractal subsets $K$  of $\R^m$. Indeed,  Mahler in the mid-eighties instigated such a line of investigation with $K$ equal to the standard middle-third Cantor set $K$. More specifically, he raised the problem:
\begin{itemize}
\item[``\!\!]{\em How close can irrational elements of Cantor's set be approximated by rational numbers (i) in Cantor's set, and (ii)  by rational numbers not in Cantor's set?''}
\end{itemize}
Mahler's problem has inspired a body of fundamental research.  In short, the solution to  (i) in its simplest form  (see \cite[Theorem~1]{MR2295506}) states:  as close as you like!  Indeed, this is a straightforward consequence of a complete Khintchine-Jarn\'{\i}k theorem (see Theorem~\ref{khijarmeas} ) in which the denominators of the rational approximates are restricted to powers of three.   Regarding (ii),  the situation  is substantially more subtle.  Building on the pioneering work of Khalil $\&$ Luethi \cite{MR4574663},   Bénard, He $\&$ Zhang \cite{bénard2025khintchinedichotomyselfsimilarmeasures}  have  recently  established a complete analogue of Khinchine's theorem for $\mu( W(\psi) \cap K ) $   where $\mu$ is the Cantor measure.  In fact,  they prove an analogue for all self-similar probability measures on the real line  --  this  includes natural measures supported on self-similar fractal sets such as the Cantor set.  For full details and background, we refer to \cite{bénard2025khintchinedichotomyselfsimilarmeasures}  and references within. The problem of determining the analogue of the Khintchine-Jarn\'{\i}k theorem for the $s$-dimensional Hausdorff measure of  $ W(\psi) \cap K  $ remains wide open.  Indeed,  as things stand, the  analogue of the  \JB{} theorem concerning the Hausdorff dimension  of  $ W(\tau) \cap K  $ is unsolved.

\medskip

\emph{Problem $:$ \ \   Let $K$ be the middle-third Cantor set.  Determine $ \dim  W(\tau) \cap K  $ for $\tau \ge 1$.  
}

\medskip

\noindent {In \cite{MR2295506}, it was speculated  that $ \dim  W(\tau) \cap K  = 2d_K/(\tau+1)   $ where $d_K:=\dim K=\log2/\log3$.  Although not disproved,  the following more likely statement 
was subsequently made in 
 \cite[Conjecture~1.1]{MR3500835}:  for  $\tau \ge 1$
$$
\dim(W(\tau)\cap K)=\max\left\{\;\frac{2}{\tau+1}-m_K,\ \ \frac{d_K}{\tau+1}\;\right\}\,,
$$
where $m_K:=1-d_K$ is the codimension of the Cantor set.  Note that this brings into play  the generic identity  \eqref{generic} for  relatively small $\tau$  (that is $1\le \tau< 1/m_K$)  and overall is in line with \eqref{PC}.  
}

\section{$M_0$ sets, Fourier dimension and normality}

 In this section the discussion revolves around three related  facts concerning the fundamental set $W(\tau)$  appearing in the \JB{} theorem: (i) it is an $M_0$ set, (ii) it is of full Fourier dimension and thus is a Salem set  and (iii)  it contains normal numbers. With this in mind, we start by setting the scene.  Throughout, $ X$ will be a bounded subset of $\R^n$  that supports a non-atomic probability measure $\mu$.  As usual,
 the Fourier transform of $\mu$ is defined by
$$ \widehat{\mu}(\xi) \, := \int   e^{-2\pi i \langle x\cdot\xi\rangle}  \, \mathrm{d}\mu (x)      \hspace{1cm} (  \xi\in\R^n )\, .  \vspace{2mm}$$
The set $X$ is called an {\em $M_0$-set} and $\mu$ a {\em Rajchman measure}  if $\widehat{\mu}(\xi) $  decays to zero at infinity.  Loosely speaking, Fourier decay reflects the ``smoothness'' or ``randomness'' of the measure. It indicates  that the measure is spread out in a way that dampens  coherent long-range correlations. In harmonic analysis and geometric measure theory, this decay is deeply tied to fractal and dimensional properties of the support of $\mu$.
 Indeed, a classical result of Frostman states:
\begin{equation}  \label{frost}
 \dim \, X  \   \ge  \  \dim_F  X  \, ,
\end{equation}
where $\dim_{F}(X) $ is  the Fourier dimension of $X$, defined as  the supremum over  all  $\beta \in [0,n]$  such that
$$
\widehat{\mu}(\xi)  =  O \big(|\xi|^{-\beta/2}\big)   \qquad \text{as}\quad |\xi|  \to \infty 
$$
for some probability measure $\mu$  supported on $X$.  The upshot of \eqref{frost} is that the decay rate of $\widehat{\mu}$  is bounded above by $\frac12 \dim X  $.  When we have equality in \eqref{frost}, the Fourier dimension is maximal and the set $X$ is called a \emph{Salem set}.

Every Borel subset of $\R^n$ with Hausdorff dimension zero is a Salem set. At the opposite extreme, $\R^n$ itself is a Salem set of dimension $n$. Less trivially, every hyper-sphere in $\R^n$ is a Salem set of dimension $(n-1)$.

The existence of Salem sets in $\R$ with arbitrary dimension $\alpha \in (0,1)$  was first established by Salem  \cite{MR43249}, using a random Cantor-type construction. This result was extended by Kahane \cite{MR212888}, who showed that for every $\alpha \in (0,n)$, there exists a Salem set in $\R^n$ of dimension $\alpha$. Kahane's approach is based on considering the image of compact subsets of [0,1] under certain stochastic processes. A common theme in the works of Salem and Kahane, as well as in many subsequent constructions of Salem sets
(see, for instance, \cite{MR1396622, MR2545245, MR3756896} and references within),
is their reliance on randomness. As a consequence, these constructions do not yield explicit examples of Salem sets; rather, they generate uncountable families of sets, almost all of which are Salem sets.  This is where the set $W(\tau)$ and the \JB{} theorem enter the narrative.

In the influential papers \cite{MR610711,MR640914},  Kaufman shows that the  naturally occurring  badly  and $\tau$--well  approximable  sets in Diophantine approximation are $M_0$ sets. In particular, for the latter case,   Kaufman constructed a probability measure $\mu$ supported on  $W(\tau)$ for any $\tau > 1$ satisfying the following  decay property:
\begin{equation} \label{kaufdec5}
\widehat{\mu}(t)=  |t|^{-\frac{1}{\tau+1}} \  o\left(\log|t| \right)   \hspace{6mm} {\rm as} \quad  |t|\to\infty \, .  \vspace{2mm}
\end{equation}
Thus, together with the \JB{} theorem, the upshot is   that  $W(\tau) $ is a Salem set for any $\tau > 1$.  We state this formally as the   \Ja-Besicovitch-Kaufman theorem.

\begin{theorem}[The \Ja-Besicovitch-Kaufman theorem] For any $\tau\ge1$, we have that  
$$
 \dim W(\tau) = \dim_F W(\tau)   =   {\frac2{\tau + 1}}\,.
$$
\end{theorem}

\noindent This result represents the first and, to the best of our knowledge, only   explicit (i.e., not random)  examples of  Salem sets in $\R$ of  arbitrary dimension $\alpha \in (0, 1)$.  Thus,  as with Hausdorff dimension of the set $W(\tau)$, the Fourier dimension of $W(\tau)$ also varies continuously with the parameter $\tau$ and thus attains every value in $(0,1)$.
Bluhm  \cite{Bl} subsequently generalised this statement to  arbitrary decreasing  functions $\psi$ -- see also \cite{MR1650442} in which a detailed account of Kaufman's work is given.
 As far as we are aware,  all known explicit examples of Salem sets of dimension other than $0, n-1$ or $n$ in $\R^n$  are based on Kaufman’s construction.
 Regarding, explicit examples in higher dimensions, it would be reasonable to suspect that simultaneous sets   $  \cA_{n}(\tau) $ for $\tau> 1/n$ and more generally $
\cA_{n,m}(\psi) $ for $ \tau >n/m$ give rise to Salem sets.  However, this is false.  Hambrook $\&$ Yu \cite{MR4621861} showed that the Fourier dimension of $\cA_{n,m}(\psi) $ is  $2n/(1 +\tau)$ when  $\tau > n/m$, thereby proving that $\cA_{n,m}(\psi) $  is not a Salem set when $\tau >n/m $  and $ mn >1 $ (recall the Hausdorff dimension is given by \eqref{BovDod}).  Over the last decade,  there has been significant activity in constructing  explicit Salem sets in higher dimensions.
Most notably, building on the two dimensional work of Hambrook \cite{MR3628226},   this effort  culminated in the elegant  work of Fraser $\&$  Hambrook~\cite{MR4548424},  who showed that  for $\tau > 1$, certain  ``algebraic'' variants of the simultaneous sets  $  \cA_{m}(\tau) $ are indeed Salem sets, with dimension $2n/(1 + \tau)$.  In short, in these constructions the approximating rationals in the  classical sets are replaced by ratios of elements from the ring of integers of a number field of degree $m$.
For details and a comprehensive overview of prior developments leading to this impressive result, see their well-written paper~\cite{MR4548424}  and the references therein.

We now turn our attention to the third and final fact outlined at the start of this section; namely, that the set $W(\tau)$, for $\tau > 1$, contains normal numbers.  Recall, that Borel (1909) proved that almost every  real number (with respect to Lebesgue measure) is normal.  However, this fundamental result gives no information regarding the existence of $\tau$--well approximable normal numbers,  since the set $W(\tau)$ is of zero-Lebesgue measure for $\tau > 1$. In the following let $X$ be a subset of the unit interval $\I$   that supports a non-atomic probability measure $\mu$ and let 
$\cA= (q_n)_{n\in \N} $  be an increasing sequence of natural numbers.  
The fundamental theorem of  Davenport, Erd\"{o}s $ \&$ LeVeque \cite{MR153656} in the theory of uniform distribution, shows that  the generic distribution properties  of  a sequence  $(q_nx)_{n\in \N} $  with  $x$ restricted to the support of $\mu$ are intimately related to the  decay rate of $\widehat{\mu}$.  In particular, when the sequence $\cal{A}$ is lacunary (i.e., there exists a  constant $ K > 1 $ such that $
  q_{n+1}  \ge K  \,  q_n $ for  $ n \in \N$), 
the theorem gives rise to the following elegant statement.

\vspace*{2mm}

\begin{thdelcorr}
Let $\mu$ be a probability  measure supported on a subset $F$ of
$ \, \I  $.  Let $\cA= (q_n)_{n\in \N} $ be a lacunary sequence of natural numbers. Let  $f: \N \to \R^+$ be a decreasing function such that
\begin{equation*} \label{lacdelsum}
\sum_{n=2}^\infty \frac{f(n)}{n \log n} \; < \; \infty  
   \quad { and  \ suppose \ that \ }  \quad  \widehat{\mu}(t)= O\left( f( |t|)  \right)  \ \  as \  \  |t|\to\infty \, .
\end{equation*}
Then  the sequence $(q_n x)_{n \in \N}$ is uniformly
distributed modulo one  for $\mu$--almost every $x \in F$.
\end{thdelcorr}

\noindent  The deduction of the proposition from the  theorem of  Davenport, Erd\"{o}s $ \&$ LeVeque is reasonably  straightforward -- the details can be found in \cite[Appendix~A]{MR4425845}.  Reinterpreting the corollary in terms of normal numbers, it implies that  $\mu$--almost every  number  in $X$ is normal. Thus, it provides a useful  mechanism for proving the existence of normal numbers in any  given zero-Lebesgue measure subset of real numbers.  Indeed, Proposition~\!DEL implies that if a given set $X$ supports a probability  measure $\mu $ such that  for some $\epsilon > 0$
\begin{equation} \label{lacdelmu2}
\widehat{\mu}(t)= O\left( (\log \log |t|)^{-(1+ \epsilon)}     \right)   \hspace{6mm} {\rm as}  \quad  |t|\to\infty \, , \vspace{2mm}
\end{equation}
then $\mu$--almost every number in $X$ is normal.  The upshot of this, together with the fact that for any $\tau > 1$, the set $W(\tau) $  supports
a probability measure $\mu$  satisfying the decay condition  \eqref{kaufdec5},   is the existence of  $\tau$--well approximable numbers that are normal.  Moreover, since $W(\tau)$  is a  Salem set,   it follows that the set of  $\tau$--well approximable normal numbers is of maximal Hausdorff dimension; namely $2/(1 + \tau)$.
The observation involving \eqref{lacdelmu2}  is  also key in many other settings, for example,  in proving  the existence of normal numbers  that  are either  (i)  badly approximable or (ii) Liouville. In the latter Liouville case, building on the work of Kaufman \cite{MR610711}, a suitable measure supported on the set $\cL$ of Liouville numbers was constructed by Bluhm  \cite{MR1657762}; see also \cite{MR1920005}.  Now focusing on the former case, the existence of the measure $\mu$,  supported on the set  
$$
{\Bad:=\{x\in\R:\inf_{q\in\N}\|qx\|>0\}}$$ 
of badly approximable numbers  {- a set of full Hausdorff dimension but Lebesgue measure zero -}  satisfying the required Fourier decay, was established by Kaufman \cite{MR610711}. More precisely, given an integer  $N \ge 2$, let  $F_N$ denote the set of real numbers in the unit interval whose continued fraction expansions have  partial quotients  bounded above by $N$. 
Clearly,  $F_N \subset \Bad$   and  for any $N\ge 3$, Kaufman constructed a probability measure $\mu$ supported on $F_N$ satisfying the  decay property:  $
\widehat{\mu}(t)= O\left( |t|^{-0.0007}  \right)  $  as $  |t|\to\infty $. 
 Kaufman's construction was subsequently refined  by Queffel\'ec $\&$ Ramar\'e \cite{MR2028020}. In particular, they showed that $F_2$ supports a probability measure with polynomial decay and they obtained faster decay rates with explicit parametrisation depending on $\dim F_N$. Unlike the situation for $W(\tau)$ ($\tau > 1$), it remains an open problem whether $F_N$ ($N \ge 2$) is a Salem set.

\medskip

     \emph{ Question $:$   \ \ For $N \ge 2$, is $F_N$ a Salem set?}

\medskip
 
\noindent  Indeed, the best known decay rate estimates currently fall well short of the required  $\frac12 \dim F_N$ for $F_N$ to attain Salem status.

We conclude this section by briefly highlighting two additional settings in which Kaufman’s work on the aforementioned Diophantine sets has had a notable impact. In the theory of multiplicative Diophantine approximation, a well-known conjecture of Littlewood states that
  $ \liminf_{q \to \infty} q \, ||q\alpha|| \, ||q\beta|| = 0 \;$
for  any pair of real numbers  $(\alpha, \beta)$. It is easily seen, via
the theory of continued fractions that the conjecture is true if
either $\alpha$ or $\beta$  are not in $\Bad$. {Also the conjecture is true if  $1$, $\alpha$, $\beta$ are linearly dependent over $\Q$. These can be regarded as trivial examples satisfying the conjecture.}  In view of this, the following  natural
question arises. Given $\alpha \in \Bad$, are there any {non-trivial}
$\beta \in \Bad$ for which the conjecture is true?  By exploiting the fact that $F_N$ is an $M_0$ set,  the following positive result is established in \cite{MR1819996}.
Given $\alpha \in \Bad$, there exists a subset $\G$ of $\Bad$ with
$\dim \G = 1$ such that for any $\beta \in \G$ the inequality
$
\, q \, ||q\alpha|| \,
||q\beta|| \, \leq \,  1 / \log q  \,
$
is satisfied for infinitely many $q \in \N$. For such pairs
$(\alpha,\beta) \in \Bad \times \Bad$, this obviously implies  the conjecture.  For subsequent developments, including the use of the ``Bohr set'' technology from additive combinatorics, see the most recent papers \cite{MR4742725, MR4834296} and references within.  Before moving on, given the context of this survey it is appropriate to mention that the analogue of the \JB{} theorem for the  $m$-dimensional multiplicative setup essentially follows the $(m-1)$-dimensional \JB{} theorem. It  states that for any $\tau \ge 1$
$$
 \dim  \cA^{\times}_{m}(\tau)   =   m-1 +  \frac{2}{\tau + 1}\, , 
$$
where $\cA^{\times}_{m}(\tau)$ denotes the set of  $x \in \I^m$ such that $
\|qx_1\|\cdots\|qx_m\|< q^{-\tau} $ for infinitely many $q \in \N$.  The details of this can be found in \cite{MR500931}.   For the  multiplicative analogue of the more subtle Khintchine-Jarn\'{\i}k theorem   (i.e., Theorem \ref{khijarmeas})  see \cite{MR3330347, MR3758214}. 
Turning to the second  setting,  we note that the set $F_N$ is a (non-linear) Cantor set.  In particular, within the language of dynamical systems,  it can be realised as the attractor of a conformal iterated function system (IFS).  As a consequence, Kaufman's work  has inspired a surge of interest in determining Fourier decay rates for dynamically defined measures, such as those supported on self-similar fractal sets. For recent progress in this direction, see, for example, \cite{baker2024fourierdecayl2flattening, MR4375453} and references within.

\vspace{4ex}

\noindent\textbf{Acknowledgements.}  We are in debt to  Besicovitch for his inspirational work that has provided us with a pretty good living -- maybe not a rich living in the monetary sense  but most certainly a  fruitful one! Also we thank Caroline Series for her encouragement, trust in us and  her  patience! We are also grateful to an anonymous reviewer for their careful reading and numerous helpful comments on an earlier version of this paper.

{\footnotesize
\bibliographystyle{abbrv}
\bibliography{Besicovitch}
}

\vspace{8ex}

{\small
\noindent VB: Department of Mathematics, University of York\\
\hspace*{4.8ex}Ian Wand Building, Deramore Lane,
York, YO10 5GH, UK\\
\hspace*{4.8ex}e-mail: victor.beresnevich@york.ac.uk

\vspace{1ex}

\noindent SV: Department of Mathematics, University
of York\\
\hspace*{4.8ex}Ian Wand Building, Deramore Lane,
York, YO10 5GH, UK\\
\hspace*{4.8ex}e-mail:  sanju.velani@york.ac.uk

}

\end{document}